\documentstyle[12pt]{article}
\begin{document}
\begin{center}
{\bf \large  Ricci curvature, minimal surfaces and sphere theorems}
\end{center}
\begin{center}
Ying Shen \\ Department of Mathematics\\
Texas A\&M
University\\College Station, TX 77843
\end{center}
\begin{center}
Shunhui Zhu\\Department of Mathematics\\
Dartmouth College\\Hanover, 
NH 03755 
\end{center}
\footnotetext[3]{The second author is supported in part by NSF grant. } 
\begin{abstract}
Using an analogue of Myers' theorem for minimal surfaces and three
dimensional topology, we prove the diameter sphere theorem for Ricci
curvature in dimension three and a corresponding eigenvalue pinching
theorem. This settles these two problems for closed manifolds with positive 
Ricci curvature since they are both false 
in dimensions greater than three.
\end{abstract}

\vspace{1cm}
\noindent {\bf \S  1. Introduction}

In this paper, we consider n-dimensional Riemannian manifolds with 
Ricci curvature $Ric \geq n-1$. Recall that under this condition, 
Myers' theorem implies that the diameter $\mbox{diam}\leq \pi$. 
Cheng \cite{Ch}, using an eigenvalue comparison method, showed that 
equality is achieved in Myers' theorem, i.e., $\mbox{diam} =\pi$, 
if and only if the manifold is isometric to the standard sphere. 
This improved an earlier result of Toponogov \cite{To} who obtained 
the same conclusion under the condition that the sectional curvature
$\mbox{sec}\geq 1$.

Associated to Toponogov's metric rigidity theorem, there is a 
topological stability theorem due to Grove-Shiohama \cite{GS} that 
says if $\mbox{sec} \geq 1$ and $\mbox{diam}>\pi/2$ then $M$ 
is homeomorphic to a sphere. Efforts were made to establish a similar
 topological stability result for Ricci curvature associated to 
Cheng's rigidity theorem, see \cite{Es}, \cite{GP}, \cite{Sh}, \cite{PZ}  
and \cite{Pe}, among others. All these results assume in addition 
some conditions on the sectional curvature. This turned out to be necessary, 
since Anderson \cite{An} and Otsu \cite{Ot} constructed examples in 
dimension $\geq 4$ that satisfy $Ric\geq n-1,\;\;\mbox{diam}\geq 
\pi -\epsilon$ for any $\epsilon$, but are not spheres. Thus the only 
dimension left is dimension three. Using the result in \cite{Zh}, 
Wu \cite{Wu} proved the three dimensional sphere theorem under an 
additional condition that the volume $\mbox{vol}>v>0$. Our first 
 theorem is to get rid of this volume condition thus settles 
the problem for Ricci curvature:\\

\noindent{\bf Theorem 1}{ \em There is a constant $\epsilon>0$ ($\approx 0.47
\pi$), such that 
if $M^3$ is a three dimensional manifold satisfying:
\[Ric \geq 2,\;\;\; \mbox{diam} > \pi -\epsilon,\]
then $M^3$ is diffeomorphic to $S^3$.}\\

We conjecture that the result is true for $\epsilon = \pi/2$.\\

There are similar results about the rigidity and stability for the 
first eigenvalues. Recall that Lichnerowicz \cite{Li} proved that $\lambda
_1\geq n$ 
if $Ric \geq n-1$, and when $\lambda _1=n$, Obata \cite{Ob} proved $M$ is
isometric 
to the standard sphere. Many authors tried to prove a stabilty result 
that says if $\lambda <n+\epsilon$, then it is a sphere. Under the 
stronger condition that $sec\geq 1$, Li-Zhong\cite{LZ}, Li-Treibergs \cite{LT}
proved a sphere theorem in lower dimensions, the problem was eventually solved
by 
Croke \cite{Cr}. The similar question for Ricci curvature again turns 
out to be false in dimension $\geq 4$ \cite{An}. It is a direct 
consequence of our theorem 1 and a result of Croke (Theorem B in 
\cite{Cr}) that the stability theorem is true for Ricci curvature 
in dimension three, namely:\\

\noindent {\bf Theorem 2} {\em There is a constant $\delta>0$ 
($\approx 1.91$),
such 
that if $M^3$ is a three dimensional manifold satisfying:
\[ Ric\geq 2,\;\;\; \lambda _1 < 3+\delta,\]
then $M$ is diffeomorphic to $S^3$.}\\

{\em Acknowledgment} We would like to thank Professor Richard
Schoen for bringing
the paper \cite{SY2} to our attention. The first author would like to thank
Professors Richard Schoen and Huai-Dong Cao for their stimulating
discussion and constant encouragement. The second author would like to
thank Professor Jeff Cheeger for his interest when this problem was
first discussed.\\

\noindent {\bf \S 2. The proof} \\

Our idea depends in essential ways on works of Hamilton \cite{Ha} and
Schoen-Yau \cite{SY1}. We adopt part of the proof to our situation with the 
help of three dimensional topology. This part was analogous to 
the argument in \cite{Zh}. Another main ingredient is a Myers' 
theorem for minimal surfaces essentially contained in \cite{SY2}. \\

\noindent {\bf Lemma 1} {\em For any $\epsilon >0$, there is a 
positive number $\tau (\epsilon)$ with 
$\lim _{\epsilon \rightarrow 0}\tau = 0$ such that if $M^n$ is a 
n-dimensional Riemannian manifold satisfying
\[Ric\geq n-1,\;\;\; \mbox{diam}\geq \pi -\epsilon,\]
then every element of $\pi _1(M, p)$ can be represented by 
a geodesic loop at p of length $\leq \tau$. Here $p$ is a 
point that realizes the diameter of $M$.}\\

\underline{\em Proof}   Let $p,q$ be two points in $M$ with
$\mbox{dist}(p,q)=\mbox{diam}$.  Let $\pi:\tilde{M}\rightarrow M$ be the
unversal covering map and 
$\tilde{p}\in \pi ^{-1}(p), \tilde{q}\in \pi ^{-1}(q)$. For any 
$[\sigma]\in \pi _1(M,p)$ with $\sigma $ a minimal geodesic loop at 
$p$, let $\tilde{\sigma}$ be a lifting of $\sigma$ with base point 
$\tilde{p}$, then $\tilde{\sigma}(1)\in \pi ^{-1}(p)$, and therefore
\[\mbox{dist}(\tilde{\sigma}(1), \tilde{q}) \geq \pi -\epsilon,\;\;\;
 \mbox{dist}(\tilde{p}, \tilde{q}) \geq \pi -\epsilon,\]
Let $\tau = \mbox{dist}(\tilde{p}, \tilde{\sigma}(1))$. Since the 
balls $B_{\tilde{p}}(\tau /2)$ and $B_{\tilde{\sigma}(1)}(\tau/2)$ 
are disjoint, one concludes
\[ \mbox{vol}(M) \geq \mbox{vol}(B_{\tilde{p}}(\tau/2))+
\mbox{vol}(B_{\tilde{\sigma}(1)}(\tau/2))+\mbox{vol}(B_{\tilde{q}}
(\pi-\epsilon-
\tau/2)).\]
By the relative volume comparison theorem,
\begin{eqnarray*}
1 &\geq & {\mbox{vol}(B_{\tilde{p}}(\tau/2)) \over 
\mbox{vol}(\tilde{M})}+{\mbox{vol}(B_{\tilde{\sigma}(1)}
(\tau/2))\over \mbox{vol}(\tilde{M})}+{\mbox{vol}(B_{\tilde{q}}
(\pi-\epsilon-\tau/2))\over \mbox{vol}(\tilde{M})}\\
& \geq & {\mbox{vol}(B^1(\tau/2)) \over \mbox{vol}(S^n)}+
{\mbox{vol}(B^1(\tau/2))\over \mbox{vol}(S^n)}+
{\mbox{vol}(B^1(\pi-\epsilon-\tau/2))\over \mbox{vol}(S^n)}
\end{eqnarray*}
where $B^1(r)$ is a ball of radius $r$ in the standard sphere $S^n$.
It follows that $\lim_{\epsilon \rightarrow 0} \tau =0$.
\hspace*{\fill}q.e.d.\\

For a simple closed curve $\Gamma\subset M$ that bounds a disk, 
we define
\[Rad(\Gamma)=\sup\{r:\;\; \Gamma \;\;
\mbox{does not bound a stable minimal disk in the $r$ neighborhood of }\;
\Gamma\}.\]\\

\noindent {\bf Lemma 2} {\em If $M^3$ is a three dimensional manifold 
with Ricci  curvature $Ric \geq 2$, then we have
\[Rad(\Gamma)\leq \frac {4\sqrt{2}}{9}\pi.\]}\\

\underline{\em Proof} This type of estimates was first introduced by
Scheon-Yau 
\cite{SY2} in the context 
of the size of black-holes by assuming a positive lower bound on the
scalar curvature. In the presence of a positive lower bound on the 
Ricci curvature, 
we can improve their estimate for $Rad(\Gamma)$, using a different
test function.  For completeness, 
we present a detailed proof in this case.

We denote the stable minimal surface under consideration by $\Sigma$
(In our application, we will actually have an area minimizing surface).
 By the second variation formula, one has
\[\int_{\Sigma}(|\nabla \psi|^2 -\psi^2 |h|^2 -\psi^2 Ric(\nu))\geq 0\]
for all $\psi$ vanishing at $\Gamma$, where $Ric(\nu)$ is the Ricci 
curvature of the ambient manifold in the direction $\nu$ which is the
unit normal vector to $\Sigma$ and $h$ is the second fundamental form.
If we shrink the surface $\Sigma$ a little bit, we know, 
by Theorem 1 in \cite{FS}, there is a positive function $g$ on the new 
minimal surface which is again denoted by $\Sigma$, 
such that
\begin{equation}
 \triangle g +|h|^2g +Ric(\nu)g=0.\nonumber
\end{equation}

Let $\rho$ be any number that is $< Rad(\Gamma)$. Fix a point 
$x$ in $\Sigma$ with $\mbox{dist}(x, \Gamma)\geq \rho$. Let $\gamma$ 
be any curve on $\Sigma$ from $x$ to $\Gamma$ and $f$ be a positive function 
on 
$\Sigma$. We consider the functional
\[I[\gamma]=\int _{\gamma} f(\gamma(s)) ds \]
where the $f$ is a function of $g$ which will be chosen later.

We minimize this functional over all such curves $\gamma$ (thus $f=1$ 
in the usual variational proof of Myers theorem). 
Let $\gamma$ be the minimizing curve with length $l\geq \rho$. 
The nonnegativity of the second variation implies that
\begin{eqnarray}
\int_{0}^{l}\{ \varphi^2 \triangle_{\Sigma} f- \varphi ^2\frac {d^2f}
{dt^2} +f \dot{\varphi}^2 -\varphi^2 f^{-1}|(\nabla f)^{\bot}|^2\}\nonumber \\
\geq \int_{0}^{l}\{f \varphi^2 R(e, \dot{c}, e, \dot{c})
+f\varphi^2 h(\dot{c},\dot{c})h(e,e)\}dt. \nonumber
\end{eqnarray}
where $\{e, \dot{c}\}$ is an orthonormal basis for the surface $\Sigma$,
$R(e, \dot{c}, e, \dot{c})$ is the sectional curvature of $M^3$ and $(\nabla
f)^{\bot}$ means the projection of $\nabla f$ onto $e$. 

We choose $f=g^k$ for some number $k$ which will be determined later, 
and use 
the fact that g satisfies equation (1) on $\Sigma$, we have
\begin{eqnarray}
 k(k-1)\int_{0}^{l}\varphi^2 g^{k-2}|\nabla g|^2-\int_{0}^{l}\varphi^2 
\frac {d^2g^k}{dt^2} +\int_{0}^{l}g^k \dot{\varphi}^2 -\int_{0}^{l}\varphi^2 
g^{-k}|(\nabla g^k)^{\bot}|^2 \nonumber\\
\geq \int _{o}^{l}\varphi^2 g^k
\{ k|h|^2 +k Ric(\nu) +
R(e, \dot{c}, e, \dot{c}) -(h(\dot{c}, \dot{c}))^2 \}.\nonumber
\end{eqnarray}

Replace $\varphi^2 g^k$ by $\varphi^2$ and integrate by parts, then one can
easily see that the last inequality implies
\begin{eqnarray}
 (\frac {1}{4}k^2 -k)\int_{0}^{l} \varphi ^2 |(\ln g)'|^2 + k \int_{0}^{l}
\varphi \varphi' (\ln g)' + \int_{0}^{l} \dot{\varphi}^2 \nonumber \\
\geq \int_{0}^{l}
\varphi^2 \{ k|h|^2 - (h(\dot{c}, \dot{c}))^2 +kRic(\nu) 
+R(e, \dot{c}, e, \dot{c})\} \nonumber
\end{eqnarray}

If we choose $0<k<4$, then by completing the square of the left-hand-side of 
the
last inequality, we get
\[ \frac {4}{4-k}\int_{0}^{l}\dot{\varphi}^2 \geq \int_{0}^{l}\varphi^2 
\{ k|h|^2 - (h(\dot{c}, \dot{c}))^2 +kRic(\nu) 
+R(e, \dot{c}, e, \dot{c})\}\]

In order to find the smallest possible upper bound for $l$ from the previous
inequality, it is not difficult to see that
one should choose $k=\frac {7}{4}$. Therefore, the last inequality becomes
\[ \frac {16}{9} \int_{0}^{l}\dot{\varphi}^2 \geq \int_{0}^{l}\varphi^2 \{
\frac{3}{4}Ric(\nu) + \frac {1}{2}S\}\]
where $S$ is the scalar curvature of the ambient manifold.

Since $Ric \geq 2$, we obtain
\[ \frac{32}{81}\int_{0}^{l} \dot{\varphi}^2 \geq \int_{0}^{l}  \varphi^2\]

Let $\phi = \sin {\pi \over l}s$, it follows that 
\[ l \leq \frac {4 \sqrt{2}}{9}\pi .\]
\hspace*{\fill}q.e.d.\\

\noindent{\bf Lemma 3} {\em Let $M\subset int(N)$ be two compact orientable 
three-manifolds with nonempty boundary. If $\pi _2(M)\rightarrow \pi _2(N)$ 
is trivial, then $\pi _1(M)$ is torsion free.}\\

This is lemma 3.5 in \cite{Zh}, see there for a proof.\\

\underline{\em Proof of Theorem 1}   
By the well-known theorem of Hamilton (\cite{Ha}), 
we only need to prove that $\pi _1(M,p)=\{e\}$ (We remark that in higher
dimensions, 
$M$ is not necessarily simply connected under the Ricci curvature condition.) 
We can also assume that $M$ is oriented. In fact, if $M$ is not orientable, 
we can use the same argument to the orientable double cover to conclude that 
$M$ has to be $RP^3$, which is orientable, thus a contradiction. In what
follows, 
we will assume all geodesic balls correspond to regular values of distance 
function 
(or its smoothing), thus are all manifolds with boundary. 
The main idea of the proof is to produce a minimal surface of sufficient 
big radius if the manifold is not simply connected. We will first prove the 
following,

\underline{\em Claim}:  For any r with $\tau < r< \mbox{diam}-2 \tau$, 
$\pi _1(B_p(r))$ is torsion free.

In fact, take $R$: $r< R<\mbox{diam}-\tau$. consider the inclusion 
$i: B_p(r)\rightarrow B_p(R)$, we only need to show that $i_*(\pi _2(B_p(r))$ 
is trivial. The claim then follows from Lemma 3.

We proceed by contradiction. Assume $i_*(\pi _2)$ is not trivial, 
by the sphere theorem in 3-dimensional topology, there is an embedded 
$S^2$ in $B_p(r)$ which is not null-homotopic in $B_p(R)$. This leads to 
three situations, each will lead to a contradiction, as we now discuss.

{\em Case 1.} $S^2$ does not seperate $B_p(R)$.

From Lemma 3.8 in \cite{He}, $B_p(R) = M_1\sharp M_2$, 
where $M_1$ is a $S^2$ bundle over $S^1$. In particular, 
$\pi _1(M_1) =Z$. But $M=M_1 \sharp M_3$ for some $M_3$. But 
$M_3 \supset B_q(\tau)$, so Lemma 1 implies that $\pi_1(M_3) \neq \{e\}$.
It follows 
from Van Kampen's theorem that $\pi _1 (M)$ is infinite, 
contradicting Myers theorem.

{\em Case 2.} $S^2$ seperates $B_p(R)$ into two components, 
both of which have non-empty intersection with $\partial B_p(R)$.

Denote $B_p(R) \setminus S^2 = M_1\cup M_2$, take 
$x_i \in M_i \cap \partial B_p(R)$. Let $\gamma$ be a curve 
connecting $x_1$ to $x_2$ in $B_p(R)$, then $\gamma$ intersects 
$S^2$ at a single point. Take another curve $\sigma$ from $x_2$ 
to $x_1$ in $M\setminus B_p(R)$, then $\sigma \circ \gamma$ is a
 closed loop in $M$, intersecting $S^2$ at a single point. Such 
a loop has infinite order in $\pi _1(M)$, again contradicting 
Myers theorem.

{\em Case 3.} $S^2$ seperates $B_p(R)$ into two connected components, 
$B_p(R)\setminus S^2 = M_1\cup M_2$, of which 
$M_1\cap \partial B_p(R) \neq \emptyset$ and 
$M_2\cap \partial B_p(R) = \emptyset$ (i.e., $\partial (M_2) =S^2$.)

It follows $B_p(R) = M_1\sharp M_2$. If $\pi _1(M_2) =\{e\}$, 
then $M_2$ is a homotopy three ball, thus the embedded $S^2$ is 
trivial in $\pi_2 (B_p(R))$, contradicting our choice of $S^2$. 
Therefore $\pi_1(M_2) \neq \{e\}$. But then $M=M_2\sharp M_3$ for some 
$M_3$. If $M_3$ is not simply connected, Van Kampen's theorem will 
again imply $\pi _1(M)$ is infinite, contradicting Myer's theorem. 
Thus $\pi_1(M_3) =\{e\}$. But $M_3 \supset B_q(\tau)$, therefore 
any loop in $B_q(\tau)$ is null-homotopic in $M_3$. This by Lemma 1 
implies $\pi _1(M, q)$ is trivial, contradicting the assumption. This completes
the proof of the claim.

Consider now the inclusions
\[ B_p(\tau)\stackrel{i}{\rightarrow} B_p(R) \stackrel {j} {\rightarrow}
 B_p(R+\delta)\]
where $R= \pi -\epsilon -2\tau$, and $\delta <\tau$.
Then
\[\pi_1(B_p(\tau))\stackrel{i_*}{\rightarrow} \pi _1(B_p(R) )
\stackrel {j_*}\rightarrow  \pi _1(B_p(R+\delta)).\]
By the claim, the middle group is torsion free, thus $i_*(\pi _1(B_p(\tau))$ 
is either trivial or has infinite order. The first case is not possible 
since it will imply $\pi _1(M)$ is trivial. Thus $i_*(\pi _1(B_p(\tau))$ 
is infinite.

Take a geodesic loop $\sigma$ at $p$ representing an element of 
$\pi _1(M, p)$ such that $i_*(\sigma)$ has infinite order. 
Since $\pi _1(M,p)$ is finite, there is an $m>0$ such that 
$m\cdot\sigma$ is trivial in $\pi _1(M,p)$. Purturbe this curve 
to have it embedded (still denoted by $m\cdot\sigma$). We then  solve
the Platau problem  in $M$ with $m\cdot\sigma$ as boundary to get an
area minimizing surface $\Sigma$ (a disc ). Since $i_*(\sigma)$ has
infinite order in $\pi _1(B_p(R))$, 
$m\cdot\sigma$ is not trivial in $B_p(R)$, thus 
$\Sigma \cap (M\setminus B_p(R))\neq \emptyset$, 
that is, there is $x\in \Sigma$ such that 
$\mbox{dist}(x,p) >R = \pi -\epsilon -2 \tau.$ 
This is not possible by Lemma 2. Direct computation
from the above inequality and that of Lemma 2 shows that $\epsilon \approx
0.47 \pi$ is enough.
\hspace*{\fill}q.e.d.

\end{document}